\documentclass[onecolumn]{autart}
\usepackage{amssymb}
\usepackage{amsmath}
\usepackage{graphicx}  
\usepackage{epstopdf}
\usepackage{hyperref}
\usepackage{overpic}
\usepackage{enumitem}

\usepackage{pgfplots}
\usetikzlibrary{intersections}
\usetikzlibrary{patterns}
\usetikzlibrary{positioning}
\usetikzlibrary{calc}
\usepgfplotslibrary{patchplots}

\newcommand{\tr}{\intercal}

\newtheorem{theorem}{Theorem}

\newtheorem{remark}{Remark}
\newtheorem{assumption}{Assumption}
\newtheorem{proposition}{Proposition}
\newtheorem{definition}{Definition}
\newtheorem{example}{Example}

\begin{document}

\begin{frontmatter}

\title{A set-theoretic generalization of dissipativity with applications in Tube MPC\thanksref{footnoteinfo}}

\thanks[footnoteinfo]{MEV and BH acknowledge financial support via ShanghaiTech U.~Grant F-0203-14-012. }

\author[ShanghaiTech]{Mario E.~Villanueva}\ead{meduardov@shanghaitech.edu.cn},
\author[ShanghaiTech]{Elena De Lazzari},
\author[UHannover]{Matthias A.~M\"uller}\ead{mueller@irt.uni-hannover.de},
\author[ShanghaiTech]{Boris Houska}\ead{borish@shanghaitech.edu.cn}
\address[ShanghaiTech]{School of Information Science and Technology, ShanghaiTech University, China.}
\address[UHannover]{Institute of Automatic Control, Leibniz University Hannover, Germany.}


\begin{abstract}
This paper introduces a framework for analyzing a general class of uncertain nonlinear discrete-time systems with given state-, control-, and disturbance constraints. In particular, we propose a set-theoretic generalization of the concept of dissipativity of systems that are affected by external disturbances. The corresponding theoretical developments build upon set based analysis methods and lay a general theoretical foundation for a rigorous stability analysis of economic tube model predictive controllers. Besides, we discuss practical prodecures for verifying set-dissipativity of constrained linear control systems with convex stage costs.
\end{abstract}

\end{frontmatter}

\section{Introduction}
\label{sec::introduction}

Dissipativity theory can be regarded as one of the most fundamental tools for 
analyzing the stability of control systems~\cite{Byrnes1991}. The 
origins of dissipativity theory can be traced to the work by Willems~\cite{Willems1971,Willems1972}, who analyzed the theoretical properties of dissipative
systems as well as formalized the 
concepts of energy supply and energy storage for general control systems.

Recent work on dissipativity theory has focused on its application to optimally
operated control systems. For example,~\cite{Angeli2012} established a link
between dissipativity of a control system and the existence of optimal steady-states. 
In~\cite{Faulwasser2018}, a thorough review of economic model predictive 
control (MPC) schemes is presented. Unlike standard tracking problems, 
economic MPC controllers are based on objective functions which are, in general, 
not positive definite. For such controllers a number of stability conditions 
are available~\cite{Angeli2012,Muller2015,Muller2016,Zanon2017}, which all rely 
on dissipativity theory.

In order to understand why one may wish to develop a generalization of
dissipativity for set-valued systems, one must be aware of of 
set-valued analysis~\cite{Aubin2009} and its importance in the development and analysis
of robust control methods~\cite{Bertsekas1972,Blanchini2009}. Among the various 
set-theoretic control methodologies, Tube model predictive control strategies have
been analyzed exhaustively during the past two decades~\cite{Langson2004,Mayne2005}. Here, the main idea is to replace trajectories by robust forward invariant tubes (RFITs), i.e., set-valued functions in the state space enclosing
all future system states, independently of the uncertainty realization. 
A great variety of methods for Tube MPC synthesis can be found in the overview 
article~\cite{Rakovic2012b}.

Notice that there is a large body of work regarding the stability of 
nominal (certainty-equivalent) MPC schemes~\cite{Chen1998,Gruene2009,Rawlings2009}.
Of course, if a parameterized version of a Tube MPC problem can be written as a 
standard MPC problem, such stability results can be applied. For example, in the 
so-called Rigid Tube MPC~\cite{Rakovic2012c,Zeilinger2014} one computes offline both the tube cross-section as well as an ancillary feedback law in order to add robustness margins to all constraints. Thus, in this case, the robust reformulation is equivalent to a nominal MPC scheme with tightened constraints and standard stability results for MPC can be applied~\cite{Rakovic2012b,Zeilinger2014}. 
For Rigid
Tube MPC schemes with economic objectives, stability results can be obtained using
tools from the field of dissipativity theory~\cite{Bayer2014,Bayer2018,Broomhead2015}.

As rigid tubes may be rather conservative, multiple strategies have
been proposed to increase the accuracy of RFITs. These include the use of homothetic~\cite{Rakovic2012,Rakovic2013} and elastic tube parameterizations~\cite{Rakovic2016}, 
which are based on polytopic sets with a constant, pre-specified number of facets. 
The use of ellipsoidal parameterizations,~\cite{Villanueva2017}, has also been proposed
for tube MPC. In general, the question of which set parameterization is the best has no unique answer, as Tube MPC formulations face an inherent tradeoff between 
computational tractability and conservatism~\cite{Rakovic2012b}. Roughly, whenever one attempts to increase the accuracy of the set representation, the computational procedures become more demanding in terms of their memory and run-time requirements~\cite{Houska2019}.

In this context, one of the main contributions of this paper is the development
of a rigorous mathematical framework for the stability analysis of a 
rather general class of set-valued control systems. Towards this aim, 
Section~\ref{sec::CostToTravel} introduces a set-based 
generalization of cost-to-travel functions, which have originally been developed for 
certainty-equivalent control systems~\cite{Houska2017}. 
Section~\ref{sec::Dissipativity} builds on this construction to propose a set-theoretic
generalization of dissipativity for a particular class of storage functions.
The practical applicability of these rather abstract
concepts is discussed in Section~\ref{sec::MPCstability}, which establishes 
set-theoretic stability conditions for a large class of Tube MPC 
controllers with possibly economic objectives ad no assumptions on the 
feedback structure. 
These controllers can be 
based, in the most general case, on parameterizations
where the set-valued cross-sections of the tube itself are free optimization 
variables. The theoretical developments of this paper are illustrated throughout
the paper using a series of academic examples. Section~\ref{sec::Conclusions} concludes the paper.

\subsection{Notation and preliminaries}
We use the symbols $\mathbb{K}^{n}$ and $\mathbb{K}^{n}_{\rm C}$ to denote the sets of compact and compact convex subsets of $\mathbb{R}^{n}$, 
respectively. The Hausdorff distance between two sets $A,B \in \mathbb K^n$ is denoted by
\begin{equation*}
d_{\rm H}(A,B) = \max \left\{ \max_{x \in A} \min_{y \in B} \Vert x-y \Vert , \,
 \max_{y \in B} \min_{x \in A} \Vert x-y \Vert \right\} .
\end{equation*}
Notice that $( \mathbb{K}^{n}, d_\mathrm{H} )$ is a metric space~\cite{Srivastava2008}.

As this paper uses functions whose arguments are sets in~$\mathbb K^{n}$, we introduce the following definitions.
\begin{definition}
\label{def::LD}
Let $\mathcal D \subseteq \mathbb{K}^{n}$ be a given domain. A function $L: \mathcal D \to \mathbb R$ is called
\begin{enumerate}
\item continuous on $\mathcal D$ if there exists for every $A \in \mathcal D$ and every $\epsilon > 0$ a $\delta > 0$ such that $|L(A) - L(B)| < \epsilon$, for all $B \in \mathcal D$ with $d_{\rm H}(A,B) \leq \delta$,
\item lower semi-continuous on $\mathcal D$ if there exists for every $A \in \mathcal D$ and every $\epsilon > 0$ a constant $\delta > 0$ such that $L(B) > L(A) - \epsilon$, for all $B \in \mathcal D$ with $d_{\rm H}(A,B) \leq \delta$, and
\item monotonous if $A \subseteq B$ implies $L(A) \leq L(B)$.
\end{enumerate}
\end{definition}

Moreover, we also introduce the generalized Hausdorff distance,
\begin{equation*}
 H( \mathcal D, \mathcal E )   =  \max \left\{ \max_{A \in \mathcal D} \min_{B \in \mathcal E} d_\mathrm{H}(A,B), \,
 \max_{B \in \mathcal E} \min_{A \in \mathcal D} d_\mathrm{H}(A,B) \right\} \; ,
\end{equation*}
which is defined for any $\mathcal D, \mathcal E \subseteq \mathbb K^n$. The
symbol $\mathcal K^{n}$ is used to denote the topological space of all nonempty 
subsets of $\mathbb K^{n}$ that are compact in $2^{\mathbb K^{n}}$---the power 
set of $\mathbb K^n$. Recall that $d_{\rm H}$ induces a metric in $\mathbb{K}^{n}$, using this 
one can show that the generalized Hausdorff distance $H$ induces a metric in $\mathcal{K}^{n}$~\cite{Rockafellar2005}. 

The following definition is useful for analyzing difference inclusions, as needed in the context of Tube MPC.
\begin{definition}
\label{def::FD}
Consider the function $F:\mathbb{K}^{n} \to \mathcal K^n$. It is called continuous
if there exists for every $\epsilon > 0$ a $\delta > 0$ such that 
\begin{equation*}
H(F(A),F(B)) < \epsilon,  
\end{equation*}
 for all $A,B \in \mathbb K^n$ with \mbox{$d_{\rm H}(A,B) \leq \delta$}.
\end{definition}

\section{Set-based cost-to-travel functions}
\label{sec::CostToTravel}
The main goal of this paper is to analyze uncertain discrete-time control systems
of the form
\begin{align}
\label{eq::system}
x_{k+1} = f( x_k, u_k, w_k ) \; .
\end{align}
Here, $x_k \in \mathbb R^{n_x}$, $u_k \in \mathbb R^{n_u}$, and $w_k \in \mathbb R^{n_w}$
denote the state, control, and disturbance vectors at time $k$.
The disturbance sequence $w$ is unknown, but assumed to take values in the given set 
$\mathbb W\in\mathbb{K}^{n_w}$. Associated state- and control constraint sets, $\mathbb X\in\mathbb{K}^{n_x}$ and $\mathbb{U}\in\mathbb{K}^{n_u}$, are also assumed to be given.

Since~\eqref{eq::system} depends on an uncertain disturbance sequence, its 
reachable set is, in general, not a singleton. Hence, Section~\ref{sec::rfit}
briefly reviews some concepts from robust forward invariance~\cite{Blanchini2009}, 
used for the analysis. 
Section~\ref{sec::ctt} introduces a novel set-theoretic generalization 
of cost-to-travel functions~\cite{Houska2017}, whose properties are analyzed in 
Section~\ref{sec::properties}.

\subsection{Difference inclusions and robust invariance}
\label{sec::rfit}

As the focus of this paper is on set-based methods for analyzing~\eqref{eq::system}, we introduce the map $F: \mathbb K^{n_x} \to \mathcal{K}^{n_x}$, given by
\begin{equation}
\label{eq::FDef}
F(A) = \left\{
B \in \mathbb K^{n_x} \, \middle | \,
\begin{aligned}
&\forall x \in A, \; \exists u \in \mathbb U: \; \forall w \in \mathbb W,\\[0.1cm]
&f(x,u,w) \in B
\end{aligned}
\right\}
\end{equation}
for all $A \in \mathbb K^{n_x}$. This transition map $F$ is the
basis for the construction of control invariant sets and tubes for~\eqref{eq::system}.

\begin{definition}\label{def::RFIT}
A sequence $X = (X_0, X_1, \ldots)$ of compact sets is called a robust forward 
invariant tube (RFIT) for~\eqref{eq::system} if it satisfies the difference inclusion
\begin{equation*}
\label{eq::RFIT}
\forall k \in \mathbb N, \quad X_{k+1} \in F(X_{k}) \; .
\end{equation*}
If $X = (X^\star, X^\star, \ldots)$ is a time-invariant RFIT, $X^\star$ is called a robust control invariant (RCI) set.
\end{definition}

Notice that $F$ maps a set to a set of sets. This notation may appear 
rather abstract on the first view, but it has the advantage that we do not have 
introduce notation for the underlying possibly set-valued feedback law and the 
associated closed-loop reachability sequences which are parametric on the 
feedback law.

\subsection{Set-based cost-to-travel functions}
\label{sec::ctt}
Let $\mathcal D \subseteq \mathbb{K}^{n}$ be a given domain and $L: \mathcal D \to \mathbb R$ a given lower semi-continuous function on $\mathcal D$. The cost-to-travel function
 $V_{\mathcal D}: {\mathcal D} \times {\mathcal D} \times \mathbb N \to 
\mathbb R \cup \{ \infty \}$
of~\eqref{eq::system} on ${\mathcal D}$ is given by
\begin{equation}
\label{eq::Vdef}
V_{{\mathcal D}}(A,B,N) =  \underset{X \in {\mathcal D}^{N+1}}{\min} \ \sum\limits_{k=0}^{N-1} L( X_k ) 
\quad \mathrm{s.t.}  \left\{
\begin{array}{l}
\forall k \in \{ 0,1, \ldots, N-1\}, \\
X_{k+1} \in F(X_{k}) \\
X_{k} \subseteq \mathbb X, \\
X_{0} = A, \; , \; X_N = B \;,
\end{array}
\right.
\end{equation}
which is defined for all sets $A,\,B \in {\mathcal D}$ and all $N \in \mathbb N$. In order to ensure that $V_{{\mathcal D}}$ is well-defined, the following assumption is needed.

\begin{assumption}
\label{ass::blanket}
We assume that the domain ${\mathcal D}$ and the functions $f$ and $L$ have the following properties:
\begin{enumerate}
\item the right-hand side function $f$ is continuous in all its arguments,
\item the set ${\mathcal D}\subseteq\mathbb{K}^{n_x}$ is closed in the metric space $(\mathbb K^{n_x}, d_H )$, and
\item the function $L: {\mathcal D} \to \mathbb{R}$ is lower semi-continuous and monotonous on ${\mathcal D}$ in the sense of Definition~\ref{def::LD}.
\end{enumerate}
\end{assumption}

\begin{proposition}
If Assumption~\ref{ass::blanket} is satisfied, then the right-hand side of~\eqref{eq::Vdef} either admits a minimizer or has an empty feasible set.
\end{proposition}

\vspace{-1.em}
\begin{pf}
First notice that $F$ is continuous in the sense of Definition~\ref{def::FD}. This is a direct consequence of the definition of $F$ in~\eqref{eq::FDef}, the 
continuity of $f$ as well as the compactness of $\mathbb U$ and
$\mathbb{W}$; see,~\cite{Aubin2009} for details. Since $\mathbb X$ is compact and ${\mathcal D}$ closed, the feasible set of~\eqref{eq::TubeMPC} is compact in $(\mathbb K^{n_x},d_{\rm H})$. Since $L$ is lower semi-continuous, the right-hand side of~\eqref{eq::Vdef} either admits a minimizer or has an empty feasible set.
\hfill\hfill\qed
\vspace{-1.em}
\end{pf}

If $A$ and $B$ are such that~\eqref{eq::Vdef} is
infeasible, we set $V_{{\mathcal D}}(A,B,N) = \infty$. This guarantees that 
the function $V_{{\mathcal D}}$ is
well-defined for all $A,B \in \mathbb K^{n_x}$.

\begin{example}
\label{ex::tutorial}
Let us consider a dynamic system given by
\begin{equation*}
f(x,u,w) = \left(
\begin{array}{c}
u \\
\frac{1}{2} x_2 + u + w
\end{array}
\right) \; ,
\end{equation*}
with $\mathbb{X} = [-5,5] \times [-5,5]$, $\mathbb U = [-5,5]$, and 
$\mathbb{W} = [-1,1]$. 
Moreover we consider the $2$-dimensional interval domain
\begin{equation*}
\mathcal D = \left\{ \; [a_1,a_2] \times [a_3,a_4] \subseteq \mathbb R^2 \; \middle | \;
\begin{aligned}
&a_1,a_2,a_3,a_4 \in \mathbb R \\
&(a_1 \leq a_2) \wedge (a_3 \leq a_4)
\end{aligned} \;
\right\}
\end{equation*}
as well as the stage cost 
\begin{equation*}
L( [a_1,a_2] \times [a_3,a_4] ) = 2a_2 
+ \frac{1}{20}\left( 3a_1^2 + a_2^2 +2a_3^2 + a_4^2 \right) \;.
\end{equation*}
In this case, the cost-to-travel function $V_{\mathcal D}(\cdot,\cdot,1)$ can be
constructed explicitly. In fact, it is given by
\begin{equation*}
\begin{aligned}
&V_{\mathcal D}(A,B,1) = 
\begin{cases}
L( [a_1,a_2] \times [a_3,a_4] ) \ &\text{if} \quad (a,b)\in G\\
\infty & \text{otherwise} \; .
\end{cases}
\end{aligned}
\end{equation*}
for all intervals $A = [a_1,a_2] \times [a_3,a_4] \in \mathcal D$ as well as all $B = [b_1,b_2] \times [b_3,b_4] \in \mathcal D$.
Here, we have used the shorthand notation
\begin{equation*}
G = \left\{ (a,b)\in\mathbb{R}^{4}
\ \middle| \
\begin{aligned}
\exists v_1&,v_2 \in [-5,5]: \\
b_3 &\leq \frac{1}{2}a_3 + v_1 - 1 \\
b_4 &\geq \frac{1}{2}a_4 + v_2 + 1 \\
a_4 &\geq 2(v_1-v_2) + a_3\\
b_1 &\leq v_1 \leq b_2 \\
b_1 &\leq v_2 \leq b_2 \\
-5 &\leq a_1 \leq  a_2 \leq 5 \\
-5 &\leq a_3 \leq  a_4 \leq 5
\end{aligned}
\right\} \;.
\end{equation*}
\end{example}

\subsection{Properties of cost-to-travel functions}
\label{sec::properties}
The following propositions summarize basic properties of the cost-to-travel
function $V_{{\mathcal D}}$.
\begin{proposition}[Monotonicity]
\label{prop::V}
If Assumption~\ref{ass::blanket} is satisfied, then
\begin{equation*}
V_{{\mathcal D}}(A,C,N) \leq V_{{\mathcal D}}(A^\prime,C,N) \quad
\text{and} \quad V_{{\mathcal D}}(A,C,N) \geq V_{{\mathcal D}}(A,C^\prime,N) 
\end{equation*}
for all sets $A,A^\prime,C,C^\prime \in {\mathcal D}$ with $A \subseteq A^\prime$ and $C \subseteq C^\prime$ and all $N \in \mathbb N$.
\end{proposition}

\begin{pf}
As discussed above, Assumption~\ref{ass::blanket} ensures that $V_{{\mathcal D}}$ is well-defined. The definition of $F$ implies that the implications
\begin{align}
C \in F(A^\prime)  \quad &\Longrightarrow \quad C \in F(A) \notag \\[0.16cm]
C \in F(A)  \quad &\Longrightarrow \quad C^\prime \in F(A) \notag
\end{align}
hold for all sets $A,A^\prime,C,C^\prime \in {\mathcal D}$ with $A \subseteq A^\prime$ and $C \subseteq C^\prime$. Moreover, Assumption~\ref{ass::blanket} requires $L$ to be monotonous; 
that is,
\begin{align}
A \subseteq A^\prime \quad \Longrightarrow \quad L(A) \subseteq L(A^\prime) \; .
\end{align}
The statement of the proposition is a direct consequence of these three implications recalling the definition of $V_{\mathcal D}$ in~\eqref{eq::Vdef}.
\hfill\hfill\qed
\end{pf}

\begin{proposition}[Continuity]
\label{prop::V2}
Let Assumption~\ref{ass::blanket} hold. The function $V_{{\mathcal D}}(\cdot,\cdot, N)$ is lower semi-continuous on its domain $$\{ (A,B) \in {\mathcal D} \times {\mathcal D} \mid V_{{\mathcal D}}(A,B, N) < \infty \} \; .$$
\end{proposition}

\begin{pf}
Assumption~\ref{ass::blanket} ensures that $F$ is continuous and
$L$ lower semi-continuous. Since $\mathbb X$ is compact, it follows, from 
standard arguments from set-valued analysis~\cite{Aubin2009}, that
$V_{{\mathcal D}}$ is lower semi-continuous. 
For example, one can use an indirect argument, as follows.

If $V_{{\mathcal D}}$ was not lower-semi-continuous, we could find a 
sequence of sets $(A_i,B_i)$ with 
\begin{equation*}
V_{{\mathcal D}}(A_i,B_i, N) < V_{{\mathcal D}}(A,B, N) - \epsilon\;,
\end{equation*}
for some $\epsilon > 0$ as well as a feasible pair $(A,B)$, such that $(A_i,B_i)$
converges to $(A,B)$ for $i \to \infty$. But this means that there exists a sequence of associated feasible points $X^i$ of~\eqref{eq::Vdef} with $A$ and $B$ replaced by $A_{i}$ and $B_{i}$; and 
$$\sum\limits_{k=0}^{N-1} L( X_k^i ) < V_{{\mathcal D}}(A,B, N) - \epsilon \; .$$
Since $\mathbb X$ is compact, this sequence must have a convergent sub-sequence, whose limit sequence $X^\infty$ is feasible too, and satisfies 
$$\sum\limits_{k=0}^{N-1} L( X_k^\infty ) \leq V_{{\mathcal D}}(A,B, N) - \epsilon \; .$$
This is a contradiction, as we have $X_0^\infty = A$ as well as $X_N^\infty = B$ by
construction. Thus, $V_{{\mathcal D}}(\cdot,\cdot, N)$ is lower semi-continuous.
\hfill\hfill\qed
\end{pf}

The set-based cost-to-travel functions $\mathcal V_{{\mathcal D}}$, satisfies the following functional equation.

\begin{proposition}[Functional equation]
\label{prop::FunctionalEquation}
Let Assumption~\ref{ass::blanket} be satisfied. Then, $V_{{\mathcal D}}$ satisfies the functional equation
\begin{align*}
& V_{{\mathcal D}}(A,C,M+N)  =  \min_{B \in {\mathcal D}} 
V_{{\mathcal D}}(A,B,M) + V_{{\mathcal D}}(B,C,N)
\end{align*}
for all $A,C \in {\mathcal D}$ and all $M,N \in \mathbb N$.
\end{proposition}

\begin{pf}
This statement follows from the definition of 
$V_{{\mathcal D}}$ and Proposition~\ref{prop::V2}. This ensures that either
a minimizer for the minimization problem over $B$ exists or that the expressions
on both sides of the functional equation are equal to~$\infty$.
\hfill\hfill\qed
\end{pf}

\section{A set-theoretic generalization of dissipativity}
\label{sec::Dissipativity}

This section introduces a generalization of dissipativity in the 
context of discrete-time set-valued inclusions.

\begin{definition}
System~\eqref{eq::system} is called set-dissipative on its domain $\mathbb X \times \mathbb U \times \mathbb W$ with respect to a given supply rate $S: {\mathcal D} \to \mathbb R$ on ${\mathcal D}$ if there exists a nonnegative storage function $\Lambda: {\mathcal D} \to \mathbb{R}_+$ such that the inequality
\begin{equation*}
\Lambda(B) - \Lambda(A) \leq S(A) \; ,
\end{equation*}
holds for all $A,B \in {\mathcal D}$ with $A,B \subseteq \mathbb X$ and $B \in F(A)$.
\end{definition}
Notice that for the special case that $\mathbb W$ is a singleton and ${\mathcal D}$
the set of singletons in $\mathbb K^{n_x}$, set-dissipativity is equivalent to 
dissipativity for deterministic systems with control-invariant supply rates, as introduced by Willems in~\cite{Willems1971,Willems1972}. To explain how set-dissipativity relates
to the ongoing developments in this paper, we introduce the following definition.

\begin{definition}
\label{def::Xstar}
A set $X^\star \in {\mathcal D}$ is called an optimal robust control invariant set if
\begin{equation*}
V_{{\mathcal D}}^\star = V_{{\mathcal D}}(X^\star,X^\star,1) = \min_{A \in {\mathcal D}} \, V_{{\mathcal D}}(A,A,1) \; .
\end{equation*}
\end{definition}

In order to ensure that $V_{{\mathcal D}}^\star$ is well-defined  the following assumption is introduced.

\begin{assumption}
\label{ass::feasibility}
The set $
\{ A \in {\mathcal D} \mid A \in F(A), \ A \subseteq \mathbb X \}$
has a non-empty interior in ${\mathcal D}$.
\end{assumption}

\begin{proposition}
\label{prop::X}
Let Assumptions~\ref{ass::blanket} and \ref{ass::feasibility} hold. Then, there exists at least one optimal robust control invariant set $X^\star \in {\mathcal D}$.
\end{proposition}

\begin{pf}
Assumption~\ref{ass::feasibility} implies that there exists at least one set $A \in {\mathcal D}$ with $A \subseteq \mathbb X$ and $A \in F(A)$, which ensures that the domain
$$\{ (A,A) \in {\mathcal D} \times {\mathcal D} \mid V_{{\mathcal D}}(A,A,1) < \infty \}$$
is non-empty. Now, the statement of this proposition is a direct consequence of Proposition~\ref{prop::V2} and Weierstrass' theorem, which can be applied here as $\mathbb X$ is compact.
\hfill\hfill\qed
\end{pf}

\begin{example}
\label{ex::tutorial2}
Consider the setting from Example~\ref{ex::tutorial}. 
Here, the optimal robust control invariant set can be found by solving 
\begin{equation}\label{eq::exampleRCI}
\min_{(a,b) \in \mathbb R^{4}} 
L([a_1,a_2]\times[a_3,a_4]) \quad \text{s.t.} \quad \left\{
\begin{aligned}
(a,b)&\in G \\
a &= b \; .
\end{aligned}
\right.
\end{equation}
Notice that~\eqref{eq::exampleRCI} is a strictly convex quadratic program with  
its unique minimizer $a^\star = b^\star = ( -1, -1, -4, 0 )^\tr$. Thus, the optimal robust control invariant set is given by the line segment $X^\star = \{ -1 \} \times [-4,0]$ 
with $V_{\mathcal D}^\star = -\frac{1}{5}$.
\end{example}

\begin{definition}
\label{def::lowerBound}
The function $V_{{\mathcal D}}(\cdot,\cdot,N)$ is called separable on ${\mathcal D}$ if it admits a non-negative separable lower bound $W: {\mathcal D} \to \mathbb R_+$ satisfying
\begin{equation*}
\forall A,B \in {\mathcal D}, \quad  
V_{{\mathcal D}}(A,B,N) - N V_{{\mathcal D}}^\star \geq W(B)-W(A) \; . \notag \\
\label{eq::lowerBound}
\end{equation*}
\end{definition}

The following lemma establishes the link between 
set-dissipativity and cost-to-travel functions.

\begin{theorem}
\label{lem::DissipativityCharacterization}
Let Assumptions~\ref{ass::blanket} and~\ref{ass::feasibility} hold. System~\eqref{eq::system} is set-dissipative on $\mathbb X \times \mathbb U \times \mathbb W$ with respect to the supply rate $S(A) = L(A)-L(X^\star)$ on ${\mathcal D}$ if and only if $V_{{\mathcal D}}(\cdot,\cdot,1)$ is separable on ${\mathcal D}$. 
\end{theorem}

\begin{pf}
Proposition~\ref{prop::X} implies that the constant offset 
$L(X^\star) = V_{{\mathcal D}}^\star < \infty$ is well-defined.
If the system~\eqref{eq::system} is set dissipative and $A$ and $B$ are such that
$V(A,B,1) < \infty$, we have
\begin{equation*}
V_{{\mathcal D}}(A,B,1) - V_{{\mathcal D}}^\star = L(A)-L(X^\star) 
\geq \Lambda(A^+) -\Lambda(A)
\end{equation*}
for all sets $A^+ \in {\mathcal D}$ with $A^+ \in F(A)$ and $A^+ \in \mathbb X$.
In particular, this inequality must hold for $A^+ = B$, which implies
\begin{equation*}
V_{{\mathcal D}}(A,B,1) - V_{{\mathcal D}}^\star \geq \Lambda(B) -\Lambda(A) \; .
\end{equation*}
This inequality also holds whenever $V(A,B,1) = \infty$. Thus,
$W = \Lambda$ is a non-negative separable lower bound of $V_{{\mathcal D}}(\cdot,\cdot,1)$
on ${\mathcal D}$. Therefore, if~\eqref{eq::system} is set-dissipative on $\mathbb X \times \mathbb U \times \mathbb W$ with respect to the supply rate $L(\cdot)-L(X^\star)$ on ${\mathcal D}$, then $V_{{\mathcal D}}(\cdot,\cdot,1)$ is separable on the domain ${\mathcal D}$.

In order to establish the converse implication, we use that $L(A) = V(A,B,1)$ for all $A,B \in \mathcal{X}$ with $A,B \subseteq \mathbb X$ and $B \in F(A)$. 
Hence, for all such $A,B$ we obtain
\begin{align}
&W(B) - W(A) \leq V_{{\mathcal D}}(A,B,1) - V_{{\mathcal D}}^\star = L(A) - L(X^\star), \nonumber 
\end{align}
which implies that~\eqref{eq::system} is set-dissipative with storage function 
$\Lambda = W$, as long as $V_{{\mathcal D}}(\cdot,\cdot,1)$ is separable on 
${\mathcal D}$ with separable lower bound $W$.
\hfill\hfill\qed
\end{pf}

\begin{example}
\label{ex::tutorial3}
Here, we continue discussing Examples~\ref{ex::tutorial} 
and~\ref{ex::tutorial2}. In this setting, the function
\begin{equation*}
W([a_1,a_2]\times[a_3,a_4]) = \left\{
\begin{array}{ll}
16 + \frac{8}{5}(a_3 - a_2) \; \; & \text{if} \; \; A \subseteq \mathbb X \\
0 & \text{otherwise}
\end{array}
\right.
\end{equation*}
happens to be a non-negative separable 
lower bound on $V_{\mathcal D}(\cdot,\cdot,1)$. Here, the offset 
$16 \geq \frac{8}{5}a_2-a_3$ is chosen such that $W$ is non-negative on $\mathbb X = [-5,5] \times [-2,2]$. To 
verify that $W$ is indeed a separable lower bound, we can compute the minimum of
the right-hand side of the inequality
\begin{equation}\label{eq::exampleW}
V_{\mathcal D}(A,B,1) - V^\star_{\mathcal D} - W(B) + W(A) \geq 0
\end{equation}
over the domain of $V_{\mathcal D}(A,B,1)$. Here, we notice that the minumum 
of the convex quadratic program
\begin{equation*}
\min_{a,b,v_1} \ L([a_1,a_2]\times[a_3,a_4]) 
- \frac{8}{5}(b_3 - b_2) + \frac{8}{5}(a_3 - a_2) \quad {\text{\rm s.t}} \left\{ 
\begin{aligned}
&b_3 \leq\frac{1}{2}a_3 + v_1 - 1\\
&v_1 \leq b_2\,,  \ a_1\leq a_2 \;,
\end{aligned}
\right.
\end{equation*}
is $-\frac{1}{5}$, with unique minimizer at
$(a^\star)^\intercal = (-1,-1,-4,0)^\intercal$, $(b^\star)^\intercal=(0,3,0,0)^\intercal$ and
$v_1^\star = 3$. Since we have $V_{\mathcal D}^\star = -\frac{1}{5}$, the 
inequality ~\eqref{eq::exampleW}
must be satisfied on the domain of $V_{\mathcal D}(\cdot,\cdot,1)$.
\end{example}

\begin{definition}
The function $V_{{\mathcal D}}(\cdot,\cdot, N)$ is called strictly separable on ${\mathcal D}$ if it is separable and the point $(X^\star,X^\star)$ is the unique minimizer of
\begin{equation*}
\min_{A,B \in {\mathcal D}} \left(V_{{\mathcal D}}(A,B,N) 
- N V_{{\mathcal D}}^\star - W(B) + W(A) \right) \; .
\end{equation*}
\end{definition}

Notice that $V_{{\mathcal D}}(\cdot,\cdot,1)$ is strictly separable if and only if~\eqref{eq::system} is set-dissipative with respect to the supply rate $S(A) = L(A) - L(X^\star)$ and the storage function $\Lambda$ is such that
\begin{equation*}
\Lambda(B) - \Lambda(A) < S(A)
\end{equation*}
for all $A,B \in {\mathcal D}$ with $A,B \subseteq \mathbb X$, $B \in F(A)$, and $(A,B) \neq (X^\star,X^\star)$.
In this sense, one may state that strict separability of $V_{{\mathcal D}}(\cdot,\cdot,1)$ is equivalent to ``strict dissipativity'' of~\eqref{eq::system}.

\section{Set-Dissipativity and Stability of Tube MPC}
\label{sec::MPCstability}

\subsection{Tube model predictive control}

Tube MPC methods proceed by solving receding-horizon optimal control problems of the form
\begin{equation}
\label{eq::TubeMPC}
\begin{alignedat}{2}
&\underset{X \in {\mathcal D}^{N+1}}{\min} \ && E(X_0) + \sum\limits_{k=0}^{N-1} L( X_k ) + M(X_N) \\[0.2cm]
&\hphantom{_X\,}\mathrm{s.t.} && \left\{
\begin{aligned}
&\forall k \in \{ 0,1, \ldots, N-1\},\\[0.1cm]
&\begin{alignedat}{2}
X_{k+1} &\in F(X_{k}) \,, \quad &&\hphantom{_N}\ z \in X_{0}, \\[0.1cm]
X_{k} &\subseteq \mathbb X \, , &&X_N \subseteq T
\end{alignedat}
\end{aligned}
\right.
\end{alignedat}
\end{equation}
with $z\in\mathbb{R}^{n_x}$ being the current state-measurement and 
$T \in {\mathcal D}$ a terminal set.
Here, $E: {\mathcal D} \to \mathbb R$, $L: {\mathcal D} \to \mathbb R$,
and $M: {\mathcal D} \to \mathbb R$ denote lower semi-continuous initial, stage, and terminal costs, respectively. 
It is well-known~\cite{Rawlings2009} that this tube MPC 
controller~\eqref{eq::TubeMPC} is recursively feasible if $T \in F(T)$ and 
$T \subseteq \mathbb X$.

\begin{remark}
If one is interested in adding a decoupled control penalty to the objective of the MPC controller,
one can always introduce discrete-time states that satisfy
\begin{equation*}
\tilde x_{k+1} = u_k \; ,
\end{equation*}
and append them to the state vector, such that the next state is equal to the current
control input. In this sense, it is not restrictive to assume
that the objective in~\eqref{eq::TubeMPC} does not explicitly depend on the control input.
\end{remark}

\begin{remark}\label{rmk::c2t-tmpc}
There is a close relation between the tube MPC problem~\eqref{eq::TubeMPC} 
and set-based cost-to-travel functions. In particular, as a direct consequence 
of Proposition~\ref{prop::FunctionalEquation},~\eqref{eq::TubeMPC} can be 
equivalently written as
\begin{equation*}
\label{eq::TubeMPCV}
\begin{alignedat}{2}
&\underset{X \in {\mathcal D}^{N+1}}{\min} \ && E(X_0) + \sum\limits_{k=0}^{N-1} V_{\mathcal D}(X_k,X_{k+1}) + M(X_N) \quad \mathrm{s.t.} && \left\{ 
\begin{aligned}
&y\in X_0 \\
&X_{N}\subseteq T\;.
\end{aligned} \right.
\end{alignedat}
\end{equation*}
\end{remark}

\subsection{Tube MPC feedback law}
Notice that, any feasible point $X$ of~\eqref{eq::TubeMPC} is 
an RFIT. Thus, we can construct a control law, 
$\mu[X]:\mathbb{N}\times\mathbb R^{n_x} \to \mathbb U$, associated to this RFIT 
such that the state of any closed-loop system
\begin{equation*}
\forall k \in \mathbb Z, \qquad  x_{k+1} = f(x_k,\mu[X](k,x_k),w_k)
\end{equation*}
satisfies the implication
\begin{equation*}
x_k \in X_k \quad \Longrightarrow \quad x_{k'} \in X_{k'}
\end{equation*}
for all $k' \geq k$ with $k, k' \in \{ 0, 1, \ldots, N \}$. This is a direct 
consequence of the definition of the transition map $F$.

\begin{remark}\label{rmk::nuk}
Consider an RFIT $X = (X_0,X_1,\ldots)$, and a point $z\in X_k$. One can evaluate
the feedback law $\mu[X](k,z)$ by solving the robust feasibility problem
\begin{equation*}
\min_{u_k} \ 0 \quad\text{s.t.}\quad f(z,u_k,w) \in X_{k+1}, \quad \forall w\in \mathbb  W
\end{equation*}
In particular, the signal $\mu[X](k,z) = u^\star_k$---with $u^\star_k$ being
a solution of the above feasibility problem, will drive $z$ to $X_{k+1}$ regardless
of the uncertainty realization. 
\end{remark}

Now, in contrast to this control law $\mu[X]$ associated to the RFIT,
the Tube MPC feedback law
$\nu: \mathbb X \to \mathbb U$ is time-invariant and given by
\begin{equation}\label{eq::tmpcnu}
\nu(z) = \mu[\Xi](0,z) \; .
\end{equation}
Here, $\Xi(z)$ denotes a minimizing sequence of~\eqref{eq::TubeMPC}
as a function of the current measurement $z$. In the following, we 
use $y = (y_0,y_1,\ldots)$ to denote the closed-loop state recursion of the Tube MPC 
controller~\eqref{eq::TubeMPC}, given by
\begin{align}
\label{eq::closedLoop}
y_{k+1} = f( y_k, \nu(y_k), w_k )
\end{align}
with $k \in \mathbb N$. That is, we set $z = y_k$, solve~\eqref{eq::TubeMPC}, 
update the system using feedback~\eqref{eq::tmpcnu} and repeat. In the next 
section we present an analysis of the stability properties of this closed-loop 
sequence using set-dissipativity.

\subsection{Stability analysis}

The goal of this section is to analyze stability of Tube MPC in the enclosure sense. Our definition of stability is motivated by the fact that the  closed-loop trajectory $y$, given by~\eqref{eq::closedLoop},
depends on the uncertainty sequence $w$.

\begin{definition}
The closed-loop state sequence $y$ is said to admit a 
stable enclosure, if there exists a sequence $Y = (Y_0,Y_1,\ldots)$ of compact sets, $Y_k \in \mathbb K^{n_x}$, such that
\begin{enumerate}
\item $y_k \in Y_k$ for all $k \in \mathbb N$, and
\item the sequence $d_{\rm H}(Y_k,X^\star)$ is stable (in the sense of Lyapunov).
\end{enumerate}
If, additionally,
\begin{equation*}
\lim_{k \to \infty} d_{\rm H}(Y_k,X^\star) =0 \; ,
\end{equation*}
then $y$ admits an asymptotically 
stable enclosure $Y$.
\end{definition}

\begin{remark}
\label{rmk::enclosure} Notice that $Y$
is not necessarily an RFIT, since the set sequence 
$Y$ is only required---under the above definition---to contain the actual 
closed-loop sequence $y$.
\end{remark}

The following theorem establishes a stability result for the Tube MPC controller~\eqref{eq::TubeMPC} under the assumption that the initial cost function $E$ is a strictly
separable lower bound of $V_{{\mathcal D}}(\cdot,\cdot, 1)$. Equivalently,
$E$ must be a storage function that establishes strict dissipativity of~\eqref{eq::system} on 
${\mathcal D}$ with respect to the supply rate $S(A) = L(A) - L(X^\star)$.
The statement is based on the additional assumption that the strictly separable lower
bounding function $E$ is also lower semi-continuous. At this point it has to be mentioned that a precise characterization of dissipative systems for which such a lower semi-continuous storage function exists, is still an open problem. However, there exist sufficient conditions under which one can assert the existence of continuous storage functions~\cite{Polushin2002}---at least for nominal (not set-valued) systems. Moreover, in the following section we will discuss a variety of cases, where one can construct continuous functions $E$ explicitly in order to arrive at a practical implementation.

\begin{theorem}
\label{thm::mainResult}
Let Assumption~\ref{ass::blanket} and~\ref{ass::feasibility} be satisfied. Let the terminal region be an optimal robust control invariant set, $T = X^\star$, and let $y_0$ be such that~\eqref{eq::TubeMPC} is feasible for $y = y_0$. If~\eqref{eq::system} is strictly set-dissipative on
$\mathbb X \times \mathbb U \times \mathbb W$ with respect to the supply rate $S(\cdot) = L(\cdot) - L(X^\star)$ on ${\mathcal D}$ with $E$ being an associated lower semi-continuous storage function and $M=0$, then the closed-loop sequence $y$ of the tube MPC controller~\eqref{eq::TubeMPC} admits an asymptotically stable enclosure.
\end{theorem}

\begin{pf}

We start the proof by constructing a 
sequence $Y = (Y_{0},Y_{1}, Y_{2},\ldots)$ as follows. 

For all $j\in\mathbb N$:\vspace{-0.75em}
\begin{enumerate}[label=(\alph*)]

\item Measure the state, $y_j$

\item\label{step::tmpc} Set $X^j = \Xi(y_j)$, where $\Xi(y_j)$ is the
optimal solution sequence of the $j$-th tube MPC problem
\begin{equation*}
\label{eq::TubeMPCAux}
\begin{alignedat}{2}
&\underset{X \in {\mathcal D}^{N+1}}{\min} \ && 
E(X^{j}_0) + \sum\limits_{k=0}^{N-1} V_{\mathcal D}
(X^{j}_k,X^{j}_{k+1},1)
\quad \mathrm{s.t.} && \left\{
\begin{aligned}
&y_j\in X^{j}_0\\
&X^{j}_N = X^\star\;.
\end{aligned} \right.
\end{alignedat}
\end{equation*}

\item Set $Y_j = X^j_0$

\item Evaluate $\nu(y_j)$, cf. Remark~\ref{rmk::nuk}, send the feedback 
signal to the system and go to (a).

\end{enumerate}
For the construction in Step~\ref{step::tmpc}, we recall the relation 
between the tube MPC problem~\eqref{eq::TubeMPC} and cost-to-travel 
functions in Remark~\ref{rmk::c2t-tmpc}.

Since we have $y_j\in X_0^j$, the relation $y_j\in Y_j$ also holds
by construction. In order to show that the sets $Y_j$ are well defined, we introduce
the shifted sequence
\begin{equation*}
\widetilde X_j 
= \left( X^j_{1}, X^{j}_{2}, \ldots, X^{j}_{N-1}, X^\star, X^\star \right) \in {\mathcal D}^{N+1}\;.
\end{equation*}
Since the inclusion $y_{j+1} \in X^{j}_{1}$ holds independently of the
uncertainty realization, $\widetilde X_j$ is a feasible point of the $(j+1)$-th 
Tube MPC problem. Thus, recursive feasibility holds and $Y_j$ is well defined.

Let $\mathcal R_{{\mathcal D}}: {\mathcal D} \times {\mathcal D} \to \mathbb R$ denote the rotated cost-to-travel function that is defined by
\begin{equation*}
\mathcal R_{{\mathcal D}}( A,B ) = E(A) - E(B) + V_{{\mathcal D}} (A,B,1) - V_{{\mathcal D}}^\star
\end{equation*}
for all $A,B \in {\mathcal D}$ such the tube MPC problem in 
Step~\ref{step::tmpc} can be written in the equivalent form
\begin{equation*}
\underset{X_j \in \mathcal D^{N+1}}{\min} \ 
\sum\limits_{k=0}^{N-1} \mathcal R_{{\mathcal D}}(X^j_{k}, 
X^j_{k+1} ) \quad 
\mathrm{s.t.} \; \; \left\{
\begin{aligned}
y_j &\in X^{j}_{0} \\[0.16cm]
X^{j}_{N} &= X^\star \; .
\end{aligned}
\right.
\end{equation*}
The key idea of this proof is to establish the claim that the function 
$\mathcal L_{\mathcal D}: {\mathcal D}^{N+1} \to \mathbb R$, given by
\begin{equation*}
\forall Z \in {\mathcal D}^{N}, \qquad \mathcal L_{\mathcal D}(Z) = \sum\limits_{k=0}^{N-1} \mathcal R_{{\mathcal D}}( Z_{k}, Z_{k+1} ) \; ,
\end{equation*}
can be used as a Lyapunov function for the iterates $X^j$ of the tube MPC 
controller. 

Our first goal is to show that the sequence $X^j$ is stable and converges to the
limit point
\begin{equation*}
\widehat X^\star = (X^\star, X^\star, \ldots, X^\star ) \in {\mathcal D}^{N+1} \; .
\end{equation*}
Let us establish the following properties of the candidate Lyapunov function 
$\mathcal L_{{\mathcal D}}$.

\begin{description}
\item[P1] The function $\mathcal L_{{\mathcal D}}$ is lower semi-continuous 
(in the sense of Definition~\ref{def::LD}). 

\item[P2] The function $\mathcal L_{{\mathcal D}}$ is positive definite, i.e.,
it satisfies $\mathcal L_{{\mathcal D}}(Z) = 0$ if and only if 
$Z = \widehat X^\star$ and $\mathcal L_{{\mathcal D}}(Z) > 0$ otherwise.
 
\item[P3] The sequence $X^j$ satisfies
\begin{equation*}
\mathcal L_{{\mathcal D}}\left(X^{j+1}\right) < 
\mathcal L_{{\mathcal D}}\left(X^j\right)
\end{equation*}
for all $j$ whenever $X^j_0 \neq X^\star$.
\end{description}

Notice that P1 follows from Proposition~\ref{prop::V2}. Moreover,~P2 follows from the definition of $\mathcal{L}_{{\mathcal D}}$ and the assumption
that $E$ is a strict separable lower bound of $V_{{\mathcal D}}(\cdot,\cdot,1)$.
Thus, it remains to establish~P3. As discussed above, the proposed tube MPC controller is recursively feasible. This implies that
\begin{align*}
\mathcal L_{{\mathcal D}}\left(X^{j+1}\right) \leq 
\mathcal L_{{\mathcal D}}\left(\widehat X^{j+1}\right) 
&= \mathcal L_{{\mathcal D}}\left(X^{j}\right) - 
\mathcal R_{{\mathcal D}}\left(X^{j}_0,X^{j}_1\right) \notag \\[0.1cm]
&< \mathcal L_{{\mathcal D}}\left(X^{j} \right)
\end{align*}
whenever $X^{j}_0 \neq X^\star$. Here, we have used our assumption that 
$\mathcal V_{{\mathcal D}}(\cdot,\cdot,1)$ is strictly dissipative, which implies
that $\mathcal R_{{\mathcal D}}\left(X^{j}_0,X^{j}_1\right)> 0$ whenever $X^{j}_0 \neq X^\star$.

These properties are sufficient to conclude that $\mathcal L_{{\mathcal D}}$
is a Lyapunov function proving asymptotic stability of $X^{j}$ to 
$\hat X^\star$ with respect to the Hausdorff metric. This implies that that the 
sequence $Y$ is an asymptotically stable enclosure of $y$, converging to 
$X^*$.\hfill\hfill\qed
\end{pf}

Similar to existing results for economic MPC schemes (see~\cite{Faulwasser2018}
and references therein) Theorem 2 establishes asymptotic
stability for the proposed Tube MPC controller under a
dissipativity condition. But---in contrast to nominal, certainty-equivalent,
economic MPC schemes---here, the storage function $E$ is not
only needed for analysis purposes. In fact, the proposed Tube
MPC controller makes explicit use of the initial cost E, as the 
initial tube is not fixed but an optimization variable.


\begin{remark}
Theorem~\ref{thm::mainResult} specializes---for simplicity of presentation---on the case $T = X^\star$ and $M(A) = 0$. However, a generalization of this stability result for any terminal region $T \in {\mathcal D}$ with $T \subseteq \mathbb X$ is possible under the additional assumption that the function $M$ is lower semi-continuous and satisfies the condition
\begin{align*}
&\forall A \subseteq T, \exists B \in F(A), \notag \\
&
\left\{
\begin{array}{l}
A,B \in {\mathcal D} \\
B \subseteq T \\
M(B)-E(B) \leq M(A) + L(A) - E(A) \; ,
\end{array}
\right.
\end{align*}
see also~\cite{Angeli2012} for details. An in-depth discussion on how to construct such set-based terminal costs is, however, beyond the scope of this paper. 
\end{remark}

\begin{example}
Let us return to the setting from Examples~\ref{ex::tutorial} 
and~\ref{ex::tutorial2}---recalling that the optimal RCI set is given by 
$X^\star = \{-1\} \times [-4,0]$. Let us attempt to set up a robust MPC 
controller without initial cost and $N = 2$, i.e. 
\begin{equation}\label{eq::tmpc1}
\min_{X \in \mathcal D^3} \ L(X_0) + L(X_1) \quad \mathrm{s.t.}\left\{ 
\begin{aligned}
&\forall k\in\{0,1\}\\
&\begin{aligned}
X_{k+1}&\in F(X_k) \\
X_k &\subseteq \mathbb X\\
z&\in X_0 \\
X_2 &= X^\star \;.
\end{aligned}
\end{aligned} \right.
\end{equation}

Using the notation established in Examples~\ref{ex::tutorial} 
and~\ref{ex::tutorial2}, the set optimization problem~\eqref{eq::tmpc1}
can be formulated as the strictly convex parametric quadratic program
\begin{equation}\label{eq::tmpc1qp}
\kern -0.8em \begin{alignedat}{2}
&\min_{a,b,c\in\mathbb{R}^{4}}  \ &&L([a_1,a_2]\times[a_3,a_4]) 
+ L([b_1,b_2]\times[b_3,b_4]) \\
&\hphantom{_{a,b}}\text{s.t} &&\left\{
\begin{alignedat}{2}
(a,b)&\in G\;, \  (b,c)&&\in G \\
c^\intercal &= x^\star\;,  \quad \ \ \ z &&\in [a_1,a_2]\times[a_3,a_4]
\end{alignedat}
\right.
\end{alignedat}
\end{equation}
with $(x^\star)^\intercal = (-1,-1,-4,0)^\intercal$. Having Remark~\ref{rmk::nuk}
in mind, we can 
introduce a decision variable $u_0\in[-5,5]$ and 
augment~\eqref{eq::tmpc1} with the constraints
\begin{equation*}
\forall w\in[-1,1], \quad f(z,u_0,w)\in[ b_1,b_2 ]\times [b_3,b_4]\;,
\end{equation*}
which hold, whenever 
\begin{equation}\label{eq::nuconst}
\begin{gathered}
b_1 \leq u_0 \leq b_2\;, \\
b_3 \leq \frac{1}{2}z_2 + u_0 -1\;, 
\quad \text{and} \quad b_4 \geq \frac{1}{2}z_2 + u_0 +1
\end{gathered}
\end{equation}
hold.

Now, the parametric optimizer of~\eqref{eq::tmpc1} (augmented 
with~\eqref{eq::nuconst}) is a piecewise linear function defined on 22 
critical regions (non-overlapping interval boxes). 
 
Let us consider the region $[-5,0]\times[-4,0]$, containing $X^\star$. 
An associated parametric optimal set sequence is given by
\begin{equation*}
\Xi_0(z) = \{ z_1\} \times [z_2, 0]  \; , \; 
\Xi_1(z) = \left\{ -\frac{1}{2}z_2 -3 \right\} \times \left[ -4 , 0 \right] \; ,
\end{equation*}
and $\Xi_2(z) = X^\star$, for all $z\in[-5,0]\times[-4,0]$. 
An optimal feedback law in this region is given by
\begin{equation*}
\forall z\in[-5,0]\times[-4,0],\quad \nu(z) = u_{0}^\star(z) = -\frac{1}{2} z_2 - 3\;.
\end{equation*}

This feedback law is recursively feasible, but unstable in the enclosure sense. 
Consider a closed-loop sequence starting at $y_{0} = (-1,-2)^\intercal$. The 
initial condition is in the optimal RCI set and  
$Y_{0} = \Xi(y_0) = \{-1\} \times [-3,0] \subset X^\star$. Now, at the next
time instance we have, by construction of the RFIT, 
$y_{1} \in \Xi_1(y_0) = \{ -2 \} \times [-4,0]$---regardless
of the uncertainty realization. Notice that $\Xi_1(y_0) \cap X^\star = \varnothing$.
Since $y_1\in Y_1$ must hold by construction, no matter how the uncertainty is
realized, the closed-loop system must be unstable in the enclosure sense. 

This instability issue can be fixed by adding the initial cost term $E = W$ from
Example~\ref{ex::tutorial3}. Now, the robust MPC formulation is given by
\begin{equation}\label{eq::tmpc2}
\min_{X \in \mathcal D^3} \ E(X_0) + \sum^{1}_{k=0}L(X_k) \quad \mathrm{s.t.}\left\{ 
\begin{aligned}
&\forall k\in\{0,1\}\\
&\begin{aligned}
X_{k+1}&\in F(X_k) \\
X_k &\subseteq \mathbb X\\
y&\in X_0 \\
X_2 &= X^\star \;.
\end{aligned}
\end{aligned} \right.
\end{equation}
Again, we can formulate this as the quadratic program
\begin{equation}\label{eq::tmpc2qp}
\kern -0.8em \begin{alignedat}{2}
&\min_{a,b,c\in\mathbb{R}^{4}} && \ W([a_1,a_2]\times[a_3,a_4]) + 
L([a_1,a_2]\times[a_3,a_4]) + L([b_1,b_2]\times[b_3,b_4]) \\[0.2cm]
&\hphantom{_{a,b}}\text{s.t} &&\left\{
\begin{alignedat}{2}
(a,b)&\in G\;, \  (b,c)&&\in G \\
c^\intercal &= x^\star\;,  \quad \ \ \ y &&\in [a_1,a_2]\times[a_3,a_4]
\end{alignedat}
\right.
\end{alignedat}
\end{equation}
augmented with the decision variable $u_0\in[-5,5]$ and the 
constraints~\eqref{eq::nuconst}. The optimizer is, again, a piecewise affine 
function defined over 24 critical regions. Figure~\ref{fig::nu0} shows
the component $u_0$ of the parametric optimizer. 

 \begin{figure}[h!]
   \centering
\vspace{1em}
 \begin{overpic}[width=0.55\columnwidth]{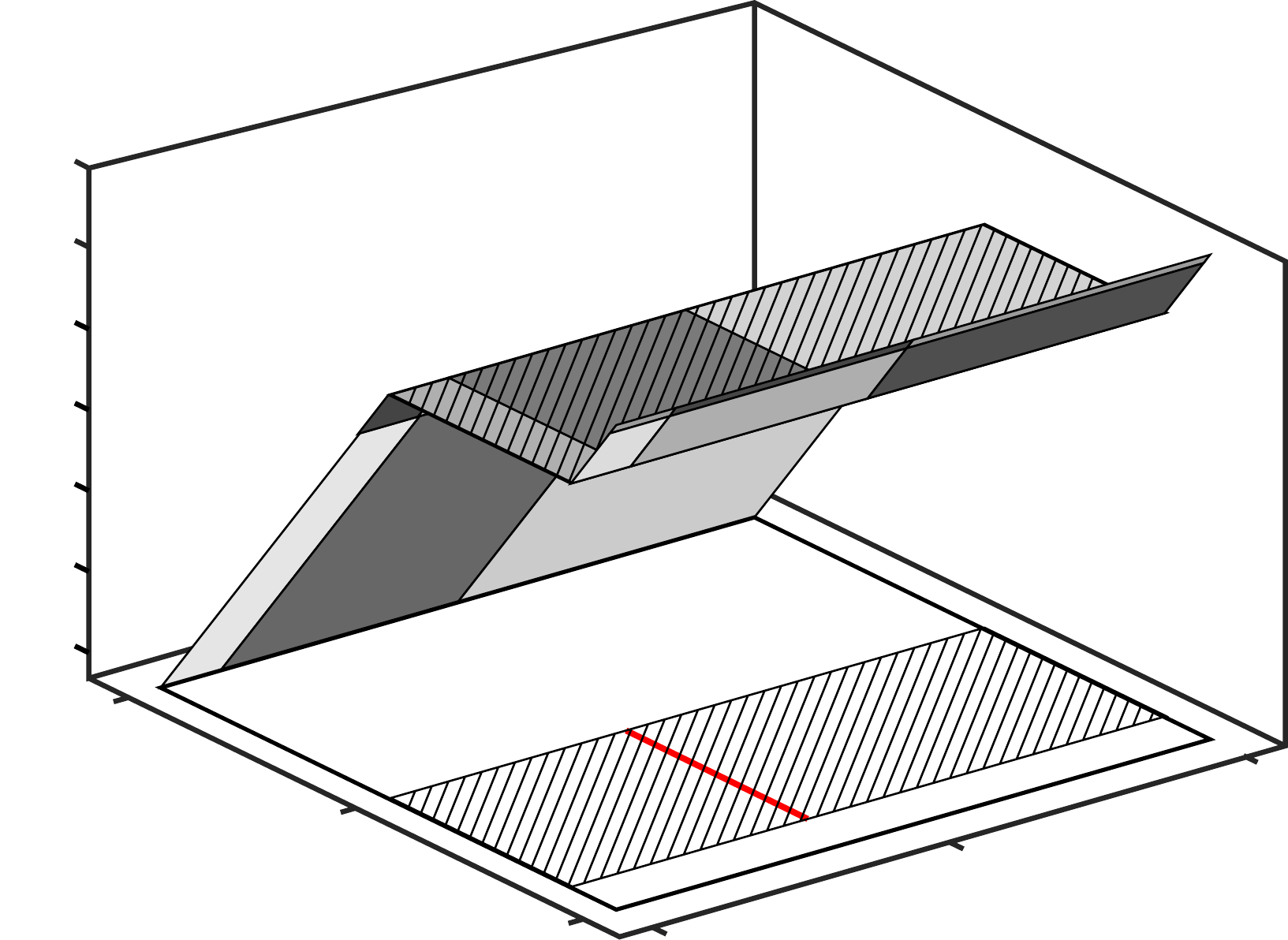}
 \put(-4,22){\scriptsize $-3.5$}
 \put(-4,28){\scriptsize $-3.0$}
 \put(-4,34){\scriptsize $-2.5$}
 \put(-4,40){\scriptsize $-2.0$}
 \put(-4,47){\scriptsize $-1.5$}
 \put(-4,53){\scriptsize $-1.0$}
 \put(-4,60){\scriptsize $-0.5$}

 \put(6,16){\scriptsize $5$}
 \put(23,7){\scriptsize $0$}
 \put(37,-1){\scriptsize $-5$}

 \put(98,10){\scriptsize $5$}
 \put(75,4){\scriptsize $0$}
 \put(49,-4){\scriptsize $-5$}

 \put(18,5){\scriptsize $z_2$}
 \put(77,0){\scriptsize $z_1$}
 \put(56,75){\scriptsize $u^\star_0(z)$}

 \end{overpic}
\vspace{1em}
 \caption{\label{fig::nu0} Component $\nu_0$ of the parametric optimizer
of~\eqref{eq::tmpc2qp}. The region $[-5,5]\times[-4,0]$ is shown hatched 
while the set $X^\star$ is shown as a red solid line.} 
 \end{figure}
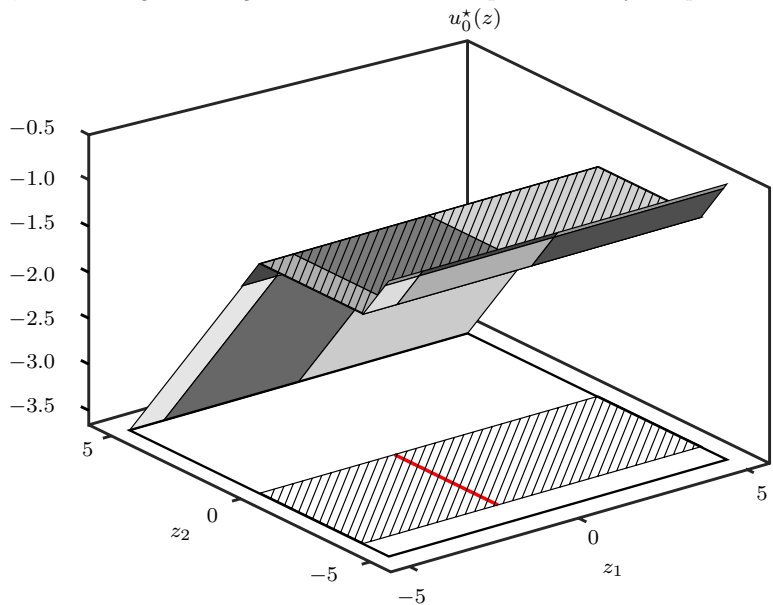

The tube MPC feedback law $\nu$, leading to the minimal stage 
cost, is given by
\begin{equation}\label{eq::tmpcnu2}
\nu(z) = u_0^\star(z) =
\begin{cases}
-\frac{1}{2} z_2 - 3 &\text{if} \ z_2\in[-5,-4]\\
-1 &\text{if} \ z_2\in[-4,0]\\
-\frac{1}{2} z_2 - 1 &\text{otherwise} \;, 
\end{cases}
\end{equation}
for all $y\in[-5,5]\times[-5,5]$. This feedback law is not only 
recursively feasible, but also asymptotically stable in the enclosure sense. 

Notice that the region $[-5,5]\times[-4,0]$---depicted with a hatched pattern in
Figure~\ref{fig::nu0}---is forward reachable in at most one step, for any initial
feasible initial condition and any $w_0\in\mathbb W$. Moreover, any 
closed-loop sequence satisfies $y_{k+1}\in X^\star$, whenever 
$y_k\in [-5,5]\times[-4,0]$---irrespective of $w_k$. Since we have 
$Y_k \subseteq X^\star$ for all $k\geq 0$, any closed loop sequence admits a 
stable enclosure. In addition, the associated optimal set sequence is given by
\begin{equation*}
\Xi_0(z) = \left\{
\begin{alignedat}{2}
[z_1,-4]&\times [-4,0] \quad &&\text{if} \ z_1\in[-5,-4] \\
\{z_1\}&\times [-4,0] &&\text{if} \ z_1\in[-4,0] \\
[0,z_1]&\times[-4,0] &&\text{otherwise} \;,
\end{alignedat}\right.
\end{equation*}
and $\Xi_1(z) = \Xi_2(z) = X^\star$  for all $z\in[-5,5]\times[0,4]$. 
Based on the previous reachability argument, it is clear that any closed
loop sequence under the feedback law~\eqref{eq::tmpcnu2} admits an asymptotically
stable enclosure. 

Figure~\ref{fig::yseq} depicts sequences closed loop sequences (blue dots with 
blue dotted lines) starting from the lower right and upper left 
corners---$(5,-5)^\intercal$ and $(-5,5)^\intercal$ respectively---of the
constraint set $\mathbb X$. The disturbance sequence has been constructed so as 
to maximize the cost. The gray sets denote the optimal sequences $\Xi(y_0)$, 
while the blue dashed lines denote the boundary of the enclosure sequences $Y$.
Notice the closed-loop system reaches the region $[-5,5]\times[-4,0]$ (hatched),
in at most 1 step, and the terminal set (red continuous
line) in at most 2 steps---remaining there, as predicted. 

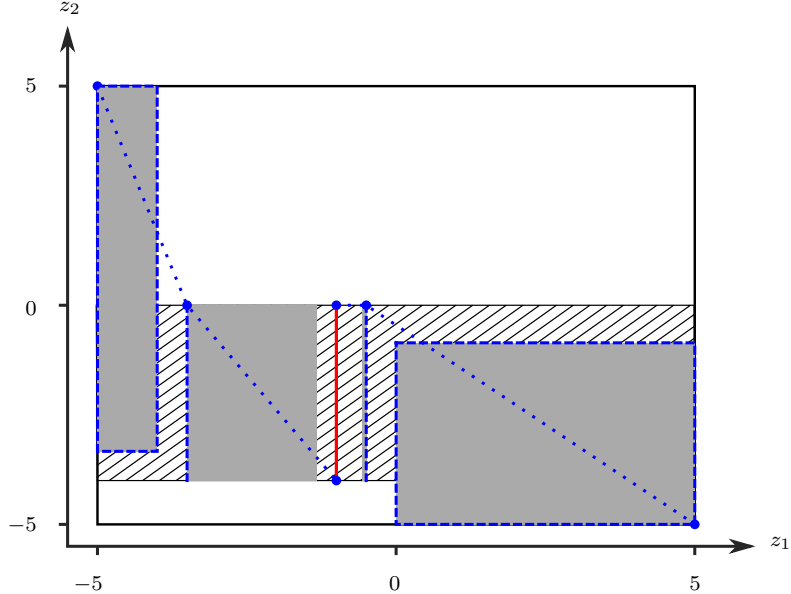
\begin{figure}
   \centering
   \vspace{0.5em}
   \begin{overpic}[width=0.55\columnwidth ]{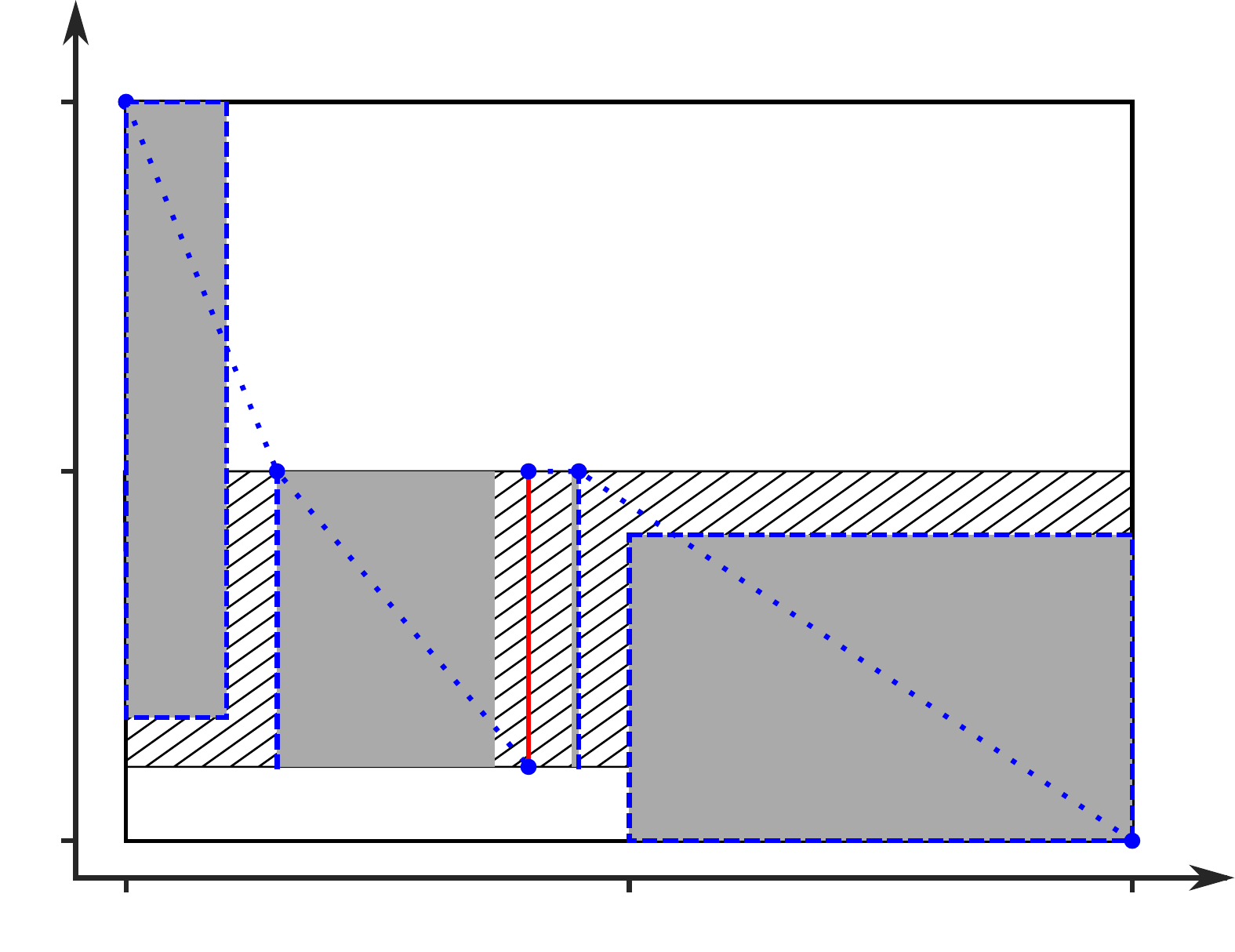}
    \put(-2,8){\scriptsize$-5$}
    \put(-2,37.5){\scriptsize$\hphantom{-}0$}
    \put(-2,68){\scriptsize$\hphantom{-}5$}

    \put(7,0){\scriptsize$-5$}
    \put(50,0){\scriptsize$0$}
    \put(91,0){\scriptsize$5$}

    \put(5,79){\scriptsize$z_2$}
    \put(102,6){\scriptsize$z_1$}
  \end{overpic}
 \caption{\label{fig::yseq} Closed-loop behavior under the feedback 
law~\eqref{eq::tmpcnu2}. The figure shows the closed-loop sequences 
$y = (y_0,y_1,y_2)$ (blue dots joined by blue dotted lines), and 
the optimal sequences $\Xi(y_0)$ (gray) starting from two initial 
conditions. The boundaries of their asymptotically stable enclosures are shown as 
blue dashed lines. The terminal set is shown as a red solid line.} 
 \end{figure}
\end{example}

\section{Conclusions}
\label{sec::Conclusions}
This paper has introduced a set theoretic generalization of dissipativity in order
to establish stability conditions for a general class of Tube MPC controllers
(cf. Theorem~\ref{thm::mainResult}). Here, the focus has been on robust MPC 
controllers, whose compact set-valued states are either entirely free optimization
variables, or belong to a finite dimensional, parametric subset ${\mathcal D}$
of all compact sets in the state space. The analysis has shown why the usual requirements for asymptotic stability of certainty-equivalent MPC controllers---namely invariance of the terminal region, a strict dissipativity condition and 
feasibility of the initial point---are not sufficient to guarantee asymptotic stability
(see the first part in Example~4). In fact, Example 4 shows that a tube MPC 
controller requires an initial cost term, which corresponds to the storage 
function in the set-dissipativity condition.

\bibliographystyle{plain}
\bibliography{references}

\end{document}